\documentstyle[11pt]{article}
\begin{document}
 \baselineskip 0.22in

\title{\textbf{New and Improved Conditions for Uniqueness of Sparsest Solutions of
Underdetermined Linear Systems}}
\author{YUN-BIN ZHAO \thanks{School of Mathematics, University of Birmingham,
Edgbaston, Birmingham B15 2TT, United Kingdom ({\tt
y.zhao.2@bham.ac.uk}). This work was partially supported by the
Engineering and Physical Sciences Research Council (EPSRC).}}

\date{(First version 22 June 2012; Revised 8 November 2012)}

\maketitle

 \textbf{Abstract.}
  The uniqueness of sparsest solutions of underdetermined linear systems  plays a
  fundamental role in the newly developed compressed sensing theory.  Several new algebraic concepts, including the sub-mutual
  coherence, scaled
mutual coherence, coherence rank, and sub-coherence rank,   are
introduced in this paper in order to develop new and improved
sufficient conditions for the uniqueness of sparsest solutions.  The
coherence rank of a matrix with normalized columns is the maximum
number of absolute entries in a row of its Gram matrix that are
equal to the mutual coherence. The main result of this paper claims
that when the coherence rank of a matrix is low,  the
mutual-coherence-based uniqueness conditions for the
 sparsest solution of a linear system can be improved. Furthermore, we prove that the  Babel-function-based
 uniqueness can be  also improved by the so-called sub-Babel function.
  Moreover, we show that the scaled-coherence-based uniqueness conditions can
  be developed, and that  the right-hand-side vector $b$ of a linear
system, the support overlap of solutions, the orthogonal matrix out
of the singular value decomposition
   of a matrix, and the range property of a transposed matrix can be
   also integrated into the criteria for the uniqueness of the sparsest solution of an underdetermined linear system.     \\



\newpage
\section{Introduction}
Consider an underdetermined system of linear equations
$$ Ax=b, $$
where $A$ is a given $m\times n$ matrix with $m <n,$  and $b\in R^m$
is a given vector. Throughout this paper, we assume that  $A$ has at
least two rows, i.e.,  $m\geq 2. $ Since the system is
underdetermined, it has infinitely many solutions. Seeking for the
sparsest solution of an underdetermined linear system has recently
become an important and common request in many applications  such
 as signal and image processing, compressed sensing, computer
vision, statistical and financial model selections,  and machine
learning (see e.g., \cite{BD06, BDE09, E10, SMF10, EK12} and the
references therein). Let $\|x\|_0$ denote the number of nonzero
components of the vector $x\in R^n.$ Then finding a sparsest
solution of a linear system amounts to the so-called
$\ell_0$-minimization problem
$$ \min \{\|x\|_0: ~ Ax=b\},$$ which is known to be NP-hard
\cite{N95, AK98}.  An intensive study of this problem has been
carried out over the past few years (see e.g., \cite{C06, D06, E10,
SMF10, EK12}),  and continues its growth in both theory and
computational methods that stimulate further cross-disciplinary
applications (see e.g., \cite{BES11, MW11, S11, EK12}). However, the
understanding of $\ell_0$-problems, from theory to computational
methods,  remains very incomplete at the moment \cite{BDE09, EK12}.
For instance, the fundamental question of when an $\ell_0$-problem
admits a unique solution has not yet addressed completely, and many
existing uniqueness claims remain restrictive so far. The main
purpose of this paper is to establish some new and improved
sufficient conditions for a linear system to have a unique sparsest
solution.

The uniqueness of sparsest solutions of underdetermined linear
systems is key to  the newly developed compressed sensing theory,
leading to a significant impact on the sparse signal and image
processing \cite{C06, D06, EK12}. So far, sufficient conditions for
the uniqueness of sparsest solutions have been developed largely via
such matrix properties as unique representation property
\cite{GR97}, spark \cite{DE03}, mutual coherence \cite{DH01},
restricted isometry property (RIP)\cite{CT05}, null space property
(NSP) \cite{CDD09, Z08}, exact recovery condition  \cite{T04}, range
property of $A^T$ \cite{Z12a}, and the verifiable conditions
\cite{JKN11, JN11}.  A crucial property for the study of uniqueness
is the spark, denoted by $\textrm{Spark}(A),$ which is the smallest
number of columns of the matrix $A$ that are linearly dependent. The
spark provides the guaranteed uniqueness of sparsest solutions, as
shown by the result below.

\vskip 0.08in  \textbf{Theorem 1.1} (\cite{DE03}).  \emph{If a
linear system $Ax=b$ has a solution $x$ satisfying $\|x\|_0
<Spark(A)/2,$ then $x$ is the unique sparsest solution to the
system.}

\vskip 0.08in

 The spark is difficult to compute. Any
computable lower bound $0< \phi(A) \leq \textrm{Spark}(A),$ however,
produces a checkable sufficient  condition for the uniqueness, such
as $\|x\|_0\leq \phi(A)/2. $ The mutual coherence of a matrix (see
the definition in section 2), denoted by $\mu(A),$  is such a
property (e.g., \cite{DH01, EB02, DE03, GN03, FN03, T04}) that
yields a computable lower bound of the spark as follows
\begin{equation} \label{bound} 1+\frac{1}{\mu(A)} \leq \textrm{Spark}
(A),\end{equation} which, together with Theorem 1.1, implies  the
following   uniqueness claim.

\vskip 0.08in  \textbf{Theorem 1.2} (\cite{DH01, F04, EB02}).
\emph{If a linear system $Ax=b$ has a solution $x$ obeying
\begin{equation} \label{coherence}\|x\|_0 <
\left(1+\frac{1}{\mu(A)}\right)/2,\end{equation}  then $x$ is the
unique sparsest solution  to the system.}

\vskip 0.08in

The condition (\ref{coherence}) is restrictive in many cases. In
\cite{T04}, the Babel function, denoted by $\mu_1(p),$ is introduced
and shown to satisfy that $ \textrm{Spark} (A) \geq \min\{p:
\mu_1(p-1) \geq 1\} \geq 1+1/\mu,$  yielding the following stronger
uniqueness condition than (\ref{coherence}).

\vskip 0.08in  \textbf{Theorem 1.3} (\cite{T04}). \emph{If a linear
system $Ax=b$ has a solution $x$ obeying
\begin{equation} \label{T-bound}\|x\|_0<
\frac{1}{2} \min\{p: \mu_1(p-1) \geq 1\},\end{equation}  then $x$ is
the unique sparsest solution to the system.}

\vskip 0.08in

Theorems 1.2 and 1.3 are valid for general matrices. When $A =
[\Phi~ \Psi]$ is a concatenation of two orthogonal matrices, Elad
and Bruckstein \cite{EB02} have shown that (\ref{coherence}) can be
improved to $\|x\|_0 < 1/\mu(A), $ and when $A$ consists of $J$
concatenated orthogonal bases, Gribonval and Nielsen \cite{GN03}
have shown that the uniqueness condition can be stated as $\|x\|_0 <
\frac{1}{2}\left(1+\frac{1}{J-1}\right)/\mu (A).$ For a general
matrix $A,$ however, it remains important, from a mathematical point
of view, to address the question: \emph{How can the bounds
(\ref{coherence}) and (\ref{T-bound})  be improved?}

 In this paper, we answer this question through the classic Brauder's Theorem.  To this end, we extract and use more properties of a
matrix than the mutual coherence, including the sub-mutual
coherence, which is the second largest inner product between two
columns of a matrix with normalized columns, and the so-called
coherence rank that turns out to be an important property for the
uniqueness of sparsest solutions. The sub-Babel function is also
taken into account in order to enhance the result of Theorem 1.3
above. One of our main results in this paper claims that \emph{for a
general matrix $A,$ when the coherence rank of $A$ is smaller than
$1/\mu(A)$, the lower bound (\ref{bound}) of Spark(A), and thus the
condition (\ref{coherence}), can be improved}.

Note that the spark of a matrix is invariant under nonsingular
scalings (see section 4 for details), but the mutual coherence is
not. This suggests a space for a further improvement of the
coherence-based uniqueness conditions. Thus we introduce the concept
of the scaled mutual coherence in section 4, which enables us to
  establish
 a  general coherence-based uniqueness condition, leading to an
  optimal lower bound of the spark in certain sense. We demonstrate
  two instant applications of the scaled mutual coherence. Note that the existing
 uniqueness  conditions  use matrix properties  only, and the role of
$b$ is completely overlooked. The sparsity of a solution, however,
can also depend on the right-hand-side vector $b$ of a linear
system. How to integrate $b$ into a uniqueness condition for
sparsest solutions  is worth  addressing (as pointed out by
Bruckstein et al. \cite{BDE09}). The first application of the scaled
mutual coherence yields such a uniqueness
 condition that depends on the property of $A$ and $b$ altogether.
 The  second  application of the scaled mutual coherence
 is a uniqueness criterion for sparsest solutions via the mutual coherence of
  an orthogonal matrix out of the singular value decomposition of $A$ (see section
  4.2 for details).

 All the above-mentioned results are developed by identifying a
  lower bound for the spark of a matrix. Any improvement of the spark condition in Theorem 1.1 leads to
  a further enhancement of these results. Although it is hard to improve Theorem 1.1 in general,
  it is possible to do so in some
situations. We show that  the
 support overlap of solutions of a linear system is the information that can  be used to
 achieve this goal (see section 5 for details). Finally,
  we introduce certain
  range properties of a matrix that guarantee
  a unique sparsest solution  to a linear system.
  Similar to  the RIP \cite{CT05, C06} and the NSP \cite{CDD09, Z08},  the range property arises naturally from the analysis of
   the uniform recovery of sparse signals (as shown in
   \cite{Z12a}).

This paper is organized as follows. We introduce several new
concepts in section 2, and use them to develop improved
  uniqueness for sparsest
solutions. The improvement of Babel-function-based condition, and
the  comparison of the existing conditions with those developed in
this paper  are given in section 3. The scaled-coherence-based
uniqueness conditions and their applications are discussed in
section 4. A further improvement of the spark condition via support
overlap of solutions is demonstrated in section 5, and the
range-property-based uniqueness is briefly introduced in section 6.

\section{Improved conditions for uniqueness of sparsest solutions}
 Let $a_i,
i=1, ..., n$ be the columns of $A.$ Recall that the mutual coherence
of $A$ (see e.g., \cite{DH01, BDE09}) is defined as
$$ \mu(A)= \max_{i\not= j} \frac{|a_i^T a_j|}{\|a_i\|_2 \cdot\|a_j\|_2}. $$
So $\mu(A) $  is the maximum absolute value of the inner product
between the normalized columns of $A.$   The lower bound
(\ref{bound}) plays a vital role in the development of the
uniqueness theory and the performance guarantee  of such algorithms
as (orthogonal) matching pursuit,  $\ell_1$-minimization, and
iterative thresholding algorithms  for the sparsest solution of
linear systems (see e.g., \cite{DH01, EB02, DE03, FN03, T04, T04b,
BDE09, E10, BD08, BD09}). Any improvement of this lower bound may
lead to an enhancement of many existing results in this field.  In
what follows, we develop an improved lower bound for spark (A) that
leads to an improved sufficient conditions for a linear system to
have a unique sparsest solution. Let us begin with a few concepts.

\subsection{Sub-mutual coherence, coherence rank, and
sub-coherence rank}

Let us  sort  the different values of the inner product $ |a_i^T
a_j|/ (\|a_i\|_2 \|a_j\|_2)$ in a descending order, and denote them
by
$$\mu^{(1)}(A) > \mu^{(2)}(A) > \cdots  > \mu^{(k)}(A). $$ Clearly,  the
largest one is $  \mu^{(1)}(A)= \mu(A),$ the mutual coherence.

\vskip 0.08in

\textbf{Definition 2.1.}  \emph{The sub-mutual coherence of $A,$
$\mu^{(2)}(A),$ is
 the second largest absolute inner
product between two normalized columns of $A:$
 $$ \mu^{(2)} (A)= \max_{i\not= j}
\left\{\frac{a_i^T a_j}{\|a_i\|_2 \cdot \|a_j\|_2}:
  ~ \frac{a_i^T a_j}{\|a_i\|_2 \cdot \|a_j\|_2} < \mu (A)  \right\}.  $$}

 In order to  introduce the next useful property of a matrix, let us consider the index set
$$ S_i(A): = \left\{j:  ~ j\not =i, ~ \frac{a_i^T a_j}{\|a_i\|_2\cdot \|a_j\|_2}
=\mu(A)  \right\}, ~i=1, ..., m. $$ Without loss of generality, we
assume  that the columns of $A$ are normalized. It is easy to see
that $S_i(A)$ counts the number of absolute entries equal to $\mu(A)
$ in $i$th row of $G=A^TA, $ the Gram matrix of $A.$  Clearly, at
least one of these sets is nonempty, since the largest absolute
entry of $G$ is equal to $\mu(A).$ Denote the cardinality of $S_i(A)
$ by $\alpha_i(A),$ i.e.,
$$\alpha_i(A) =|S_i(A)|, ~i=1, ..., m.$$ Clearly, $0\leq \alpha_i(A)
\leq n-1.$  Let
\begin{equation} \label{alpha}  \alpha (A) = \max_{1\leq i\leq m}
\alpha_i (A) = \max_{1\leq i\leq m} |S_i (A)| , \end{equation} which
is a positive number. Let $i_0$ be an index such that $$ \alpha(A)=
\alpha_{i_0} (A) =|S_{i_0}(A)|,$$ i.e., the $i_0$th row of $G$ has
the maximal number of absolute entries equal to $\mu(A).$ Then we
may define
\begin{equation}\label{beta} \beta  (A) =\max_{1\leq i\leq m, ~i\not= i_0}
\alpha_i(A) = \max_{1\leq i\leq m, ~i\not= i_0} |S_i(A)|,
\end{equation} which is the second largest number among $\alpha_i
(A), i=1, ..., m.$

\vskip 0.08in

 \textbf{Definition 2.2.}  \emph{$\alpha (A), $ given by (\ref{alpha}),  is called the coherence rank of
 $A,$ and $\beta(A),$ given by (\ref{beta}), is called sub-coherence rank of $A.$}

\vskip 0.08in

For a given matrix $A$ with normalized columns, both $\alpha (A)$
and $\beta(A)$  can be easily obtained through its Gram matrix
$G=A^TA $ or its absolute Gram matrix, denoted by $\textrm{abs}(G).$
By  the definition of $\mu(A),$ there exists at least one
off-diagonal absolute entry of $A$, say $|G_{ij}|$ (in $i$th row),
which is equal to $\mu(A). $ By the symmetry of $G$, we also have
$|G_{ji}|=\mu(A)$ (in $j$th row of $G$). Thus the symmetry of $G$
implies that $\beta(A)\geq 1.$  So, for any matrix $A,$  we have the
relation
\begin{equation} \label{aaa}   1\leq \beta (A) \leq \alpha
(A).\end{equation}

 Geometrically,  $\alpha(A)$ can be called the
\emph{Equiangle} of $A $  in the sense that it is the maximum number
of columns of $A$ that have
 the same largest angle with respect to a column, say the $i_0$-column, of $A.$

\vskip 0.08in

 \textbf{Remark 2.3.}   When all  columns of $A$ are
generated by a single vector,  then $\mu(A)=1$ and $\alpha(A) =
\beta (A) = n-1.$   When $A$ has at least two independent columns
(not all columns are generated by a single vector), then $\alpha(A)
< n-1. $ For the concatenation of two orthogonal bases $A=[\Psi ~
\Phi],$ where $ \Psi, \Phi$ are $m\times m$ orthogonal matrices, we
see that $\alpha (A) \leq n/2=m. $  As we have pointed out, all
$\mu(A), \mu^{(2)}(A), \alpha(A) $ and $\beta(A) $ can be obtained
straightaway from the Gram matrix of $A. $  For example, when $A$ is
given by \begin{equation} \label{Remark} A=\left[
       \begin{array}{cccccc}
   -0.9802  &  0.1  &  0.3521  &  0.9239  &  0.9239   &  0.7405
   \\
   -1.8282     &    0  &  1.0365 &    0.3827  & -0.3827  &
   -1.6821\\
    0.3269    &     0  &  1.3563   &      0  &       0 &  -0.2949
\end{array}
     \right], \end{equation}
     then the Gram matrix of the normalized  $A$
 $$
 G = \left[\begin{array}{cccccc}
    1   &  -0.4668 &   -0.4908 &  -0.7644 &  -0.0981  &  0.5763
    \\
   -0.4668  &  1   &  0.2020 &   0.9239  &  0.9239  &  0.3978 \\
   -0.4908  &  0.2020  & 1  &  0.4142 &  -0.0409  &  -0.5803\\
   -0.7644  &  0.9239  &  0.4142 &   1 &   0.7071 &   0.0217\\
   -0.0981  &  0.9239  & -0.0409  &  0.7071  &  1  &  0.7134\\
    0.5763  &  0.3978  & -0.5803  &  0.0217  &  0.7134  &  1
             \end{array}
           \right], $$ from which
we see that  $ \mu(A)=0.9239 >  \mu^{(2)}(A)= 0.7644 , $ and $
\alpha(A) = 2> \beta(A) =1.$

\subsection{Coherence-rank-based  lower bounds for Spark(A)}

Let us first recall the Brauer's theorem \cite{B46} (see also
Theorem 2.3 in \cite{V00}), concerning the estimate of eigenvalues
of a matrix. Let $\sigma(A):=\{\lambda: \lambda\textrm{ is an
eigenvalue of }A\}$ be the spectrum of $A.$

\vskip 0.08in

\textbf{Theorem 2.4} (Brauer \cite{B46}). \emph{Let $A=(a_{ij})$ be
an $N\times N$ matrix with $N\geq 2.$ Then, if $\lambda$ is an
eigenvalue of $A$, there is a pair $(r, q)$ of positive integers
with $r\not= q $ $(1\leq r,q\leq N)$ such that
$$ |\lambda-a_{rr}|\cdot |\lambda-a_{qq}| \leq \Delta_r
\Delta_q,   \textrm{ where } \Delta_i:= \sum_{j=1, j\not =i
}^N|a_{ij}| \textrm{ for }1\leq i\leq N.$$ Hence if $ K_{ij}(A)=\{z:
|z-a_{ii}|\cdot |z-a_{jj}| \leq \Delta_i\Delta_j\}$ for $i\not=j,$
then $ \sigma (A)\subseteq \bigcup_{i\not =j}^N K_{ij}(A).$ }

\vskip 0.08in

We make use of this classic theorem to prove the following result,
which turns out to be an improved version of (\ref{bound}) when the
coherence rank is low.

\vskip 0.08in

\textbf{Theorem 2.5.} \emph{Let $A\in R^{m\times n}$ be a matrix
with $m< n,$ and let $\alpha (A) $ and $\beta(A)$ be defined by
(\ref{alpha}) and (\ref{beta}), respectively. Suppose that one of
the following conditions holds:  (i) $\alpha (A)  <
\frac{1}{\mu(A)}; $  (ii) $\alpha (A)  \leq \frac{1}{\mu(A)} $ and
$\beta(A) < \alpha(A).$
 Then $\mu^{(2)}(A) >0 $ and
\begin{equation}\label{lowerbound} \textrm{Spark}(A) \geq 1 +  \frac{ 2\left[1-\alpha(A) \beta(A) \widetilde{\mu}(A)^2 \right]}
 {\mu^{(2)}(A) \left\{ \widetilde{\mu} (A)(\alpha(A)+\beta(A))
    + \sqrt{
 \left[\widetilde{\mu}(A) (\alpha(A)-\beta(A))
 \right]^2 +4}  \right\} } , \end{equation}
 where $\widetilde{\mu} (A) : =  \mu (A)- \mu^{(2)} (A).$}

\vskip 0.08in

 \emph{Proof.} Normalizing the columns of a matrix does
not affect any of the $\textrm{Spark}(A),$  $\mu(A), $
$\mu^{(2)}(A),$ $\alpha(A)$ and $\beta(A).$   Thus, without loss of
generality, we assume that all columns of  $A$ have unit
$\ell_2$-norms. Let $p=\textrm{Spark}(A).$ By the definition of
spark, there exist $p$ columns of $A$ that are linearly dependent.
Let $A_S$ be the submatrix consisting of these $p$ columns. Without
lost of generality, we assume $A_S=(a_1,a_2, ..., a_p).$ Thus the
$p\times p$ matrix $G_{SS}= A_S^T A_S$ is singular, since the
columns of $A_S$ are linearly dependent. Note that all diagonal
entries of $G_{SS}$ are equal to 1, and all off-diagonal absolute
entries are less than or equal to $\mu(A). $ Under either condition
(i) or (ii), we have $\alpha (A) \leq \frac{1}{\mu(A)}. $ Hence it
follows from (\ref{bound}) that
$$ 1+\alpha(A) \leq  1+\frac{1}{\mu(A)} \leq \textrm{Spark} (A). $$
So, $\alpha (A) \leq \textrm{Spark}(A)-1 =p-1.$ Note that $G_{SS}$
is a $p\times p $ matrix. Thus in every row of $G_{SS},$ there exist
at most $\alpha(A)$ absolute entries  equal to $\mu(A), $ and the
remaining $(p-1)-\alpha (A) $ absolute entries are less than or
equal to $\mu^{(2)} (A) . $ By the singularity of $G_{SS}$, $\lambda
=0$ is an eigenvalue of $G_{SS}.$ Note that the entries of $G_{SS}$
are given by $G_{ij}= a_i^Ta_j$ where $ i, j $ $=1,..., p.$ Thus by
Theorem 2.4, there  exist two different rows, say $i$th and $j$th
rows $(i\not= j)$, such that
\begin{equation} \label{relation-1}|0-G_{ii}|\cdot |0-G_{jj}| \leq
\Delta_i \Delta_j = \left(\sum_{k=1, k\not= i}^p  |a_i^T a_k|
\right)
 \left(\sum_{k=1, k\not= j}^p |a_j^T a_k|\right),  \end{equation}
 where $G_{ii}= G_{jj} =1$ are two diagonal entries of $G_{SS}.$
By the definition of $\alpha(A) $ and $ \beta (A) $, one of these
two rows contains at most $\alpha(A)$ entries with absolute values
equal to $\mu(A),$ and the next row
 contains at most $\beta(A)$ entries with absolute values equal to
 $\mu(A).$ The remaining entries in these rows are less than or equal
 to $\mu^{(2)}(A).$ Therefore,
 \begin{eqnarray} \label{relation-2}
 \left(\sum_{k=1,k\not= i}^p |a_i^T
a_k| \right)
 \left(\sum_{k=1, k\not= j}^p |a_j^T a_k|\right)
 &   \leq  &   \left[ \alpha(A) \mu(A) +(p-1-\alpha(A)) \mu^{(2)} (A)
 \right]
    \cdot
 {\Big [} \beta (A)  \mu(A)  \nonumber\\
    & &  +(p-1-\beta(A)) \mu^{(2)} (A)  {\Big ]}.
  \end{eqnarray}
    Combining  (\ref{relation-1}) and (\ref{relation-2}) leads to
  \begin{eqnarray*}  1  & \leq &   \left[ \alpha(A) \mu(A) + (p-1-\alpha(A)) \mu^{(2)} (A)
  \right]\cdot
 \left[ \beta (A)  \mu(A)
     +(p-1-\beta(A))  \mu^{(2)} (A)  \right] \\
  & = & \left[ \alpha(A) \widetilde{\mu}(A)  + (p-1) \mu^{(2)} (A)
  \right] \cdot
    \left[ \beta (A)   \widetilde{\mu}(A)
     +(p-1 )  \mu^{(2)} (A) \right],
     \end{eqnarray*}  where $\widetilde{\mu}(A) := \mu(A)-\mu^{(2)}(A). $  By rearranging
     terms, the inequality above can be written as
\begin{equation} \label{quadratic}  \left[(p-1)\mu^{(2)} (A)\right]^2+ \left[(p-1) \mu^{(2)}
(A)\right]  (\alpha(A)+\beta(A)) \widetilde{\mu}(A)  + \alpha(A)
\beta(A) \widetilde{\mu}(A)^2 -1 \geq 0.
\end{equation}
We now show that $ \mu^{(2)} (A) \not = 0.$ In fact, if  $ \mu^{(2)}
(A) =0, $ then (\ref{quadratic})  is reduced to $ \alpha(A) \beta(A)
\mu(A)^2 \geq 1,$ which contradicts both conditions (i) and (ii). In
fact,  each of conditions (i) and (ii) implies that $\alpha(A)
\beta(A) \mu(A)^2 < 1.$ Thus $ \mu^{(2)} (A)$ is positive. Note that
the quadratic equation (in $t$)
$$ t^2+ t    (\alpha(A)+\beta(A)) \widetilde{\mu}(A)  +
\alpha(A) \beta(A) \widetilde{\mu}(A)^2 -1 =0
$$ has only one positive root. So  it follows from (\ref{quadratic}) that
\begin{eqnarray*}  & & (p-1)\mu^{(2)}(A)   \\
 &   & \geq    \frac{ - ( \alpha(A)
+\beta(A))  \widetilde{\mu}(A) + \sqrt{\left[
\widetilde{\mu}(A)(\alpha(A) +\beta(A))\right]^2-4 (\alpha(A)
\beta(A) \widetilde{\mu}(A)^2 -1)
  }  } {2} \\
 & &  =   \frac{ - ( \alpha(A)
+\beta(A))  \widetilde{\mu}(A) + \sqrt{\left[(\alpha(A) -\beta(A))
 \widetilde{\mu}(A) \right]^2 +4 }  } {2} \\
&  & =   \frac{ 2\left[1-\alpha(A) \beta(A) \widetilde{\mu}(A)^2
\right]  }  { (\alpha(A)+\beta(A))
  \widetilde{\mu}(A)    + \sqrt{
 \left[(\alpha(A)-\beta(A))
 \widetilde{\mu}(A)\right]^2 +4} },
 \end{eqnarray*}
 which is exactly the relation (\ref{lowerbound}).  ~~ $\Box$

\vskip 0.08in

 The next proposition shows that (\ref{lowerbound}) is an improved lower bound
 for $\textrm{Spark }(A)$ under the condition  of Theorem
 2.5.

 \vskip 0.08in

 \textbf{Proposition 2.6.} \emph{Let $ \Psi \left(\alpha (A), \beta(A), \mu(A), \mu^{(2)}(A)\right)$ denote
  the right-hand side of the inequality (\ref{lowerbound}).  When  $\alpha (A) <
\frac{1}{\mu(A)},$ we have $$  \Psi \left(\alpha (A), \beta(A),
\mu(A), \mu^{(2)}(A)\right)   \geq
   \left(1+\frac{1}{\mu(A)}\right) +
  \left(\frac{1}{\mu^{(2)}(A)} -\frac{1}{\mu(A)} \right) (1-\alpha(A) \mu(A)). $$
When $\alpha (A) \leq \frac{1}{\mu(A)}$ and $\beta(A)<\alpha(A),$ we
have \begin{eqnarray*}  \Psi \left(\alpha (A), \beta(A), \mu(A),
\mu^{(2)}(A)\right)   & \geq &
   \left(1+\frac{1}{\mu(A)}\right) +
  \left(\frac{1}{\mu^{(2)}(A)} -\frac{1}{\mu(A)} \right) (1-\alpha(A) \mu(A)) \\
   &  &  + \frac{
 \alpha(A)
\widetilde{\mu}(A)^2  }
 {\mu^{(2)}(A) ( 1+\alpha(A)\widetilde{\mu} (A)  ) } ,
  \end{eqnarray*} }
where $\widetilde{\mu}(A) =\mu(A)-\mu^{(2)}(A).$

\vskip 0.08in

 \emph{Proof.} By using the fact $\sqrt{a^2+b^2} \leq a + b$ for any $a, b\geq
 0,$ we have
\begin{eqnarray} \label{TTT}  & &  \Psi \left(\alpha (A), \beta(A), \mu(A), \mu^{(2)}(A)\right) -1 \nonumber \\
 & &  =  \frac{ 2\left(1-\alpha(A) \beta(A)
\widetilde{\mu}(A)^2 \right)}
 {\mu^{(2)}(A) \left\{ \widetilde{\mu} (A)(\alpha(A)+\beta(A))
    + \sqrt{
 \left[\widetilde{\mu}(A) (\alpha(A)-\beta(A))
 \right]^2 +4}  \right\} }  \nonumber \\
 &   &  \geq   \frac{ 2\left(1-\alpha(A) \beta(A)
\widetilde{\mu}(A)^2 \right) }
 {\mu^{(2)}(A) \left\{ \widetilde{\mu} (A)(\alpha(A)+\beta(A))
    +
 \left[\widetilde{\mu}(A) (\alpha(A)-\beta(A))
 \right] +2  \right\} }  \nonumber \\
 & & = \frac{ 1-\alpha(A) \beta(A)
\widetilde{\mu}(A)^2  }
 {\mu^{(2)}(A) ( 1+\alpha(A)\widetilde{\mu} (A)  ) }.
 \end{eqnarray}

Case 1: $\alpha (A) < \frac{1}{\mu(A)}. $ In this case, by
(\ref{aaa}), i.e., $\beta(A) \leq \alpha(A),$ it follows from
(\ref{TTT}) that
 \begin{eqnarray*}  \Psi \left(\alpha (A), \beta(A), \mu(A), \mu^{(2)}(A)\right)
 -1
  & \geq &  \frac{ 1-\alpha(A)^2
\widetilde{\mu}(A)^2  }
 {\mu^{(2)}(A) ( 1+\alpha(A)\widetilde{\mu} (A)  ) } \\
 & = & \frac{1- \alpha(A) \widetilde{\mu}(A) }{\mu^{(2)}(A)} \\
 &    = &
\frac{1}{\mu(A)} +
  \left(\frac{1}{\mu^{(2)}(A)} -\frac{1}{\mu(A)} \right) (1-\alpha(A)
  \mu(A)).
  \end{eqnarray*}

Case 2: $\alpha (A) \leq  \frac{1}{\mu(A)} $ and $\beta(A) < \alpha
(A). $  In this case, since $\beta(A)\leq \alpha(A)-1,$ it follows
again from (\ref{TTT}) that
\begin{eqnarray*}  \Psi \left(\alpha (A), \beta(A), \mu(A), \mu^{(2)}(A)\right)
 -1
  & \geq &  \frac{ 1-\alpha(A) (\alpha(A)-1)
\widetilde{\mu}(A)^2  }
 {\mu^{(2)}(A) ( 1+\alpha(A)\widetilde{\mu} (A)  ) } \\
 & = & \frac{1- \alpha(A) \widetilde{\mu}(A) }{\mu^{(2)}(A)} + \frac{
 \alpha(A)
\widetilde{\mu}(A)^2  }
 {\mu^{(2)}(A) ( 1+\alpha(A)\widetilde{\mu} (A)  ) } \\
 &    = &
\frac{1}{\mu(A)} +
  \left(\frac{1}{\mu^{(2)}(A)} -\frac{1}{\mu(A)} \right) (1-\alpha(A)
  \mu(A)) \\
  & &  + \frac{
 \alpha(A)
\widetilde{\mu}(A)^2  }
 {\mu^{(2)}(A) ( 1+\alpha(A)\widetilde{\mu} (A)  ) },
  \end{eqnarray*}
as desired. ~~ $\Box$

Under the first case above, we see that
$$
   \left(\frac{1}{\mu^{(2)}(A)} -\frac{1}{\mu(A)} \right) (1-\alpha(A) \mu(A)
  ) > 0 , $$
and under the second case, we have
$$
  \left(\frac{1}{\mu^{(2)}(A)} -\frac{1}{\mu(A)} \right) (1-\alpha(A) \mu(A)) \\
      + \frac{
 \alpha(A)
\widetilde{\mu}(A)^2  }
 {\mu^{(2)}(A) ( 1+\alpha(A)\widetilde{\mu} (A)  ) } >0 . $$

\vskip 0.08in

Thus the next corollary follows immediately from  Proposition 2.6.

\vskip 0.08in

 \textbf{Corollary 2.7.}  \emph{Under the condition
of Theorem 2.5, we have }
$$ \Psi \left(\alpha (A), \beta(A), \mu(A), \mu^{(2)}(A)\right)  >
1+\frac{1}{\mu(A)}. $$

Therefore, the lower bound of spark given by (\ref{lowerbound}) does
 improve the bound  (\ref{bound}) when the coherence rank, $\alpha(A),$ is
 small. Proposition 2.6 also indicates explicitly  how much this
 improvement can be made at least.

 If the Gram matrix
$G$ of the normalized $A$ has two rows containing $\alpha(A)$
entries with absolute values equal to $\mu(A),$  then $\alpha(A)
=\beta(A),$ in which case  the lower bound (\ref{lowerbound}) can be
  simplified to $$  \Psi \left(\alpha (A), \beta(A),
\mu(A), \mu^{(2)}(A)\right) =
   \left(1+\frac{1}{\mu(A)}\right) +
  \left(\frac{1}{\mu^{(2)}(A)} -\frac{1}{\mu(A)} \right) (1-\alpha(A) \mu(A)). $$
Note that  $G$ has at most one absolute entry equal to $\mu(A)$ in
its every row if and only if
  $\alpha (A) = \beta(A) =1. $ In this special case, the condition $\alpha(A)<
1/\mu(A)$ holds trivially when $\mu(A)<1.$
 Thus, the next corollary follows immediately from  Theorem 2.5.

\vskip 0.08in

 \textbf{Corollary 2.8.}  \emph{Let $A\in R^{m\times n} $ be a matrix with  $m<n.$ If $\mu(A) <1 $ and $\alpha(A) =1,$  then $\mu^{(2)}(A)
>0,$ and
$$ \textrm{Spark}(A)   \geq      1+\frac{1}{\mu(A)} +
  \left(\frac{1}{\mu^{(2)}(A)} -\frac{1}{\mu(A)} \right) (1- \mu(A)
  ).$$
}

Although Corollary 2.8 deals with a special case from a mathematical
point of view, many matrices satisfy the property $\alpha(A)=1$
together with $\mu(A)<1.$ Numerical experiments show that  when a
matrix is randomly generated, the coherence rank of the matrix is
most likely
 to be 1. In fact,  the case $\alpha(A)\geq 2$ arises only when
$A$ has at least two columns, each of which has
 the same angle  to  a column of the matrix, and
such an angle is the largest one between a pair of columns of $A.$
 This phenomenon  indicates that
the coherence rank of a matrix is usually low in practice, typically
$\alpha(A)=1. $

\subsection{Uniqueness via   coherence and  coherence rank}

 Consider the class of
matrices
\begin{eqnarray} \label{MMM} {\cal M} &  = & \left\{A\in R^{m\times n}:  \textrm{ either }  \alpha (A)  \leq
\frac{1}{\mu(A)}\textrm{ and } \beta(A) < \alpha(A),  \textrm{ or }
\alpha (A)  < \frac{1}{\mu(A)} \right\} \nonumber\\
 & =  &  {\cal M}_1 \cup  {\cal M}_2,
\end{eqnarray}
where $$ {\cal M}_1 = \left\{A\in R^{m\times n}:   ~ \alpha (A)  <
\frac{1}{\mu(A)} \right\}, ~  {\cal M}_2 = \left\{A\in R^{m\times
n}: ~\alpha (A)  \leq \frac{1}{\mu(A)}\textrm{ and  }\beta(A) <
\alpha(A)
  \right\}.$$

  We now state
  the  main uniqueness claim of this section.

  \vskip 0.08in

\textbf{Theorem 2.9. } \emph{Let $A \in {\cal M}, $ defined by
(\ref{MMM}). If the system $Ax=b$ has a solution $x$ obeying
\begin{equation}\label{unique-a}
\|x\|_0 < \frac{1}{2} \left[1 + \frac{ 2\left(1-\alpha(A) \beta(A)
\widetilde{\mu}(A)^2 \right)}
 {\mu^{(2)}(A) \left\{ \widetilde{\mu} (A)(\alpha(A)+\beta(A))
    + \sqrt{
 \left[\widetilde{\mu}(A) (\alpha(A)-\beta(A))
 \right]^2 +4}  \right\} }\right], \end{equation}
 where $\widetilde{\mu} (A) :=  \mu (A)- \mu^{(2)} (A),$ then $x$ is
 the unique sparsest solution to the linear system.}

 \vskip 0.08in

This result follows instantly from Theorems 2.5  and 1.1.  As shown
by Proposition 2.6, condition (\ref{unique-a}) has improved the
well-known condition (\ref{coherence}) when $A$ is in class ${\cal
M}.$ This improvement is achieved by using the sub-mutual coherence
$\mu^{(2)}(A)$ together with (sub-)coherence rank, instead of
$\mu(A)$ only. Note that $\alpha(A), \beta(A), \mu(A)$ and
$\mu^{(2)}(A)$  can be obtained straightforward from the Gram matrix
$G= A^TA.  $ Thus the bound (\ref{unique-a}) can be easily computed.

\vskip 0.08in

By Theorem 2.5 and Proposition 2.6, we  obtain the next result.

\vskip 0.08in

 \textbf{Theorem 2.10.}  (i) \emph{Let $A \in {\cal M}_1,$ defined by (\ref{MMM}).
If the  system $Ax=b$ has a solution $x$ obeying
\begin{eqnarray}\label{unique-b}
\|x\|_0   <     \frac{1}{2} \left[ 1+\frac{1}{\mu(A)} +
  \left(\frac{1}{\mu^{(2)}(A)} -\frac{1}{\mu(A)} \right) (1-\alpha(A)
  \mu(A))
  \right]  ,
\end{eqnarray} then $x$ is the unique sparsest solution of the
linear system.}

(ii) \emph{Let $A \in {\cal M}_2,$ defined by (\ref{MMM}). If the
system $Ax=b$ has a solution $x$ obeying
\begin{eqnarray}\label{ZZZZ}
\|x\|_0   <      \frac{1}{2} \left[ 1+\frac{1}{\mu(A)} +
  \left(\frac{1}{\mu^{(2)}(A)} -\frac{1}{\mu(A)} \right) (1-\alpha(A)
  \mu(A)) +\frac{
 \alpha(A)
\widetilde{\mu}(A)^2  }
 {\mu^{(2)}(A) ( 1+\alpha(A)\widetilde{\mu} (A)  ) }
  \right]  ,
\end{eqnarray} then $x$ is the unique sparsest solution of the
linear system.}

(iii)  \emph{Let $A$ be a matrix with $\mu(A) <1$ and
$\alpha(A)=1.$
 Then the solution of $Ax=b$ obeying
\begin{equation} \label{Coherence-rank-1} \|x\|_0    <    \frac{1}{2} \left[ 1+\frac{1}{\mu(A)} +
  \left(\frac{1}{\mu^{(2)}(A)} -\frac{1}{\mu(A)} \right) (1-
  \mu(A))
  \right]
\end{equation} is the unique sparsest solution of the linear
system.}

\vskip 0.08in

Result (iii) of the above theorem  shows that for coherence-rank-1
matrices, the uniqueness criterion (\ref{coherence}) can be always
improved to (\ref{Coherence-rank-1}). As we have pointed out,
matrices (especially the randomly generated ones) are largely
coherence-rank-1, unless the matrix is particularly designed.

\vskip 0.08in

\textbf{Example 2.11.}
 Consider a randomly generated  $A$ below and the absolute Gram matrix  of its column-normalized
 counterpart
 $$ A=\left[
      \begin{array}{rrrr}
        0.0010 & -0.7998 & -0.6002 & 0.0717 \\
        0.8001 & -0.3558 & 0.4798 & -0.1913 \\
        0.5999 & 0.4801 & -0.6398 & -0.6412 \\
      \end{array}
    \right], ~ \textrm{abs}(G)=  \left[
                \begin{array}{cccc}
                  1 & 0.0025 & 0.0005 & 0.7989 \\
                  0.0025 & 1 & 0.0022 & 0.4422 \\
                  0.005 & 0.0022 & 1 & 0.4093 \\
                  0.7989 & 0.4422 & 0.4093 & 1 \\
                \end{array}
              \right].
 $$ From  $\textrm{abs}(G),$  we see that $\alpha(A) =\beta(A) = 1, $ $ \mu(A)
 =0.7989,$ and
 $
 \mu^{(2)}(A) = 0.4422.$ Note that  $\textrm{Spark}(A)/2 =2$ for this
 example.
 The standard mutual bound (\ref{coherence})
 is $(1+\frac{1}{\mu(A)})/2=1.1258, $ which is improved to
 1.2274 by (\ref{Coherence-rank-1}).

 \section{Improvement of Babel-function-based uniqueness}

 Let $A\in R^{m\times n}$ be a matrix with normalized columns. Tropp
 \cite{T04} introduced the so-called Babel-function defined as
 $$\mu_1(q) =\max_{\Lambda, |\Lambda|=q} \max_{j\not\in \Lambda}
 \sum_{i\in \Lambda} |a_i^T a_j|$$  where $a_k, k=1,...,n, $ are the columns of
 $A,$ and $\Lambda$ is some subset of $\{1,..., n\}.$  By this function, the following lower bound for spark is  obtained (see \cite{T04}):
 \begin{equation} \label{Tropp}  \textrm{Spark}(A) \geq \min_{1\leq q\leq n} \{ q:  \mu_1(q-1)\geq 1\}.
 \end{equation}
The Babel function can be equivalently defined/computed in terms of
the Gram matrix $G=A^T A. $ In fact, sorting every row of
$\textrm{abs}(G)$ in descending order yields the matrix $\widehat{G}
=(\widehat{G}_{ij}) $ with the first column equal to the vector of
ones, consisting of the diagonal entries of $G.$  Therefore, as
pointed out in \cite{E10}, the Babel function can be written as
\begin{equation} \label{index}
\mu_1(q)= \max_{1\leq k\leq m} \sum_{j=2}^{q+1} |\widehat{G}_{kj}| =
\sum_{j=2}^{q+1} |\widehat{G}_{k_0 j}|, \end{equation}  where $k_0$
denotes an index such that the above maximum is achieved. Since
$\mu_1(q-1) \leq (q-1)\mu(A),$ it is evident that
$$ \min_{1\leq q\leq n} \{ q:  \mu_1(q-1)\geq 1\} \geq
1+\frac{1}{\mu(A)}.$$ So the lower bond given by (\ref{Tropp}) is an
enhanced version of (\ref{bound}). Some immediate questions arise:
\emph{Can we compare the  lower bounds (\ref{Tropp}) and
(\ref{lowerbound})?
 Can
the lower bounds (\ref{Tropp}) and (\ref{lowerbound}) be further
improved?}

We first address the second question above, by showing that
 the Babel-function-based bound (\ref{Tropp}) can be further improved by using
  the so-called sub-Babel function. Again,
Brauer's Theorem plays a fundamental role in deriving such an
enhanced result. The  sub-Babel function, denoted by
$\mu_1^{(2)}(q)$, is defined as \begin{equation} \label{sub-Babel}
\mu_1^{(2)} (q)= \max_{1\leq k\leq m, k\not= k_0} \sum_{j=2}^{q+1}
|\widehat{G}_{kj}|,\end{equation} where $k_0$ is determined in
(\ref{index}). Clearly, we have
\begin{equation}\label{inequality}  \mu_1^{(2)} (q) \leq \mu_1(q) ~~ \textrm{ for any }1\leq q \leq
n-1.\end{equation} We have the following improved version  of
(\ref{Tropp}).

\vskip 0.08in

\textbf{Theorem 3.1.} \emph{ For any matrix $A\in R^{m\times n},$ we have
   \begin{equation} \label{Zhao}  \textrm{Spark}(A) \geq \min_{1\leq q\leq n}
\left\{ q:  ~~ \mu_1(q-1) \cdot \mu_1^{(2)} (q-1) \geq 1\right\}.
 \end{equation}}

\emph{ Proof.} Let $p=\textrm{Spark}(A).$ Then there exist $p$
columns of $A$ that are linearly dependent. Without lost of
generality, we assume $A_S=(a_1,a_2, ..., a_p)$ is the submatrix
consisting of these $p$ columns.  Since the columns of $A_S$ are
linearly dependent and normalized,  the $p\times p$ matrix $G_{SS}=
A_S^T A_S$ is singular, and  all diagonal entries of $G_{SS}$ are
equal to 1. Thus by Theorem 2.4 (Brauer's Theorem), for any
eigenvalue $\lambda$ of $G_{SS},$ there exist two different rows,
say $i$th and $j$th rows $(i\not= j)$, such that
\begin{equation} \label{rrrr}|\lambda -G_{ii}|\cdot |\lambda-G_{jj}| \leq
\Delta_i \Delta_j = \left(\sum_{k=1, k\not= i}^p  |a_i^T a_k|
\right)
 \left(\sum_{k=1, k\not= j}^p |a_j^T a_k|\right),  \end{equation}
 where $G_{ii}= G_{jj} =1$ are two diagonal entries of $G_{SS}.$
By the definition of Babel and sub-Babel functions, we see that
$$ \max\{\Delta_i,  \Delta_j\} \leq \mu_1(p-1),  ~~\min \{\Delta_i,  \Delta_j\} \leq \mu_1^{(2)} (p-1). $$
Thus  it follows from (\ref{rrrr}) that
$$  (\lambda-1)^2 \leq \Delta_i \Delta_j = \max\{\Delta_i,
\Delta_j\}  \cdot \min\{\Delta_i,  \Delta_j\}  \leq \mu_1(p-1)
\cdot\mu_1^{(2)} (p-1). $$ In particular, since $\lambda=0$ is an
eigenvalue of $G_{SS}$, we have \begin{equation} \label{IIII}
\mu_1(p-1)\cdot \mu_1^{(2)} (p-1) \geq 1. \end{equation}  So
  $p=\textrm{Spark}(A)$ implies that $p$ must satisfy (\ref{IIII}).
  Therefore,
$$ \textrm{Spark}(A) =p \geq \min_{1\leq q\leq n} \left\{q: ~~ \mu_1(q-1)
\cdot\mu_1^{(2)} (q-1) \geq 1\right\},$$ as desired. ~~ $\Box.$

\vskip 0.08in

 The next proposition shows that the lower bound (\ref{Zhao}) is an improved version of (\ref{Tropp}).

\vskip 0.08in

 \textbf{Proposition 3.2.}  Denote by
$$q^*=  \min_{1\leq q\leq n}  \left\{q: ~~
\mu_1(q-1) \cdot\mu_1^{(2)} (q-1) \geq 1\right\}, ~ \widehat{q}  =
\min_{1\leq q\leq n} \left\{q: ~~ \mu_1(q-1)   \geq 1\right\}. $$
Then $q^* \geq \widehat{q}.$ In particular, if $
\mu_1^{(2)}(\widehat{q}-1) < \frac{1}{\mu_1(\widehat{q}-1)},$ then
$q^*> \widehat{q}.$

\vskip 0.08in

 \emph{Proof.} By the definition of $q^*,$ we see that
   $ \mu_1(q^*-1) \cdot\mu_1^{(2)} (q^*-1)
\geq 1 . $ This, together with (\ref{inequality}), implies that
$\mu_1(q^*-1) \geq 1. $ Thus $$q^* \geq  \min_{1\leq q\leq n} \{ q:
~ \mu_1(q-1)  \geq 1\} =\widehat{q}. $$ We now further show that
this inequality holds strictly when the value of the sub-Babel
function are relatively small in the sense that $
\mu_1^{(2)}(\widehat{q}-1) < \frac{1}{\mu_1(\widehat{q}-1)}. $ In
fact, under this condition,
          we have
$$ \mu_1 (\widehat{q}-1)\cdot \mu_1^{(2)} (\widehat{q}-1) < 1.  $$ Note that both $ \mu_1 (q-1)$ and
$ \mu_1^{(2)} (q-1)$ are increasing functions in $q.$ The inequality
above shows that when $ \mu_1(q-1)\cdot \mu_1^{(2)} (q-1) \geq 1$,
we must have $q>\widehat{q}. $ Therefore,
$$ q^* =\min_{1\leq i\leq n} \{q: ~\mu_1(q-1)\cdot \mu_1^{(2)} (q-1) \geq 1\} > \widehat{q},$$
which shows that (\ref{Zhao}) improves (\ref{Tropp}) for this case.
~~ $\Box.$

\vskip 0.08in

The next proposition indicates that when the coherence rank of $A$
is relatively small, bound (\ref{Zhao}) is also an improved version
of (\ref{lowerbound}).

\vskip 0.08in

\textbf{Proposition 3.3.} Let $A \in R^{m\times n} $ be a given
matrix. Let $q^*$ be defined as in Proposition 3.2. If $\alpha(A) <
1/\mu(A)$ and $\alpha(A) \leq q^*-1,$ then
$$q^* \geq 1 + \frac{
2\left[1-\alpha(A) \beta(A) \widetilde{\mu}(A)^2 \right]}
 {\mu^{(2)}(A) \left\{ \widetilde{\mu} (A)(\alpha(A)+\beta(A))
    + \sqrt{
 \left[\widetilde{\mu}(A) (\alpha(A)-\beta(A))
 \right]^2 +4}  \right\} },$$
 where $\widetilde{\mu} (A) : =  \mu (A)- \mu^{(2)} (A). $

\vskip 0.08in

\emph{Proof.}  Since $\alpha (A) \leq q^*-1, $ by the definition of
$\alpha (A)$ and $\beta(A),$ it follows from (\ref{index}) and
(\ref{sub-Babel}) that
 $$ \mu_1(q^*-1) \leq \alpha(A) \mu(A) + (q^*-1-\alpha(A))
\mu^{(2)} (A), $$
$$ \mu_1^{(2)}
(q^*-1) \leq \beta (A)  \mu(A)
     +(q^*-1-\beta(A))  \mu^{(2)} (A).$$
These relations, together with  the definition of $q^*,$  imply that
 \begin{eqnarray*} 1
& \leq  &  \mu_1(q^*-1) \cdot  \mu_1^{(2)}
(q^*-1) \\
 & \leq  & \left[ \alpha(A) \mu(A) + (q^*-1-\alpha(A))
\mu^{(2)} (A)
  \right]\cdot
 \left[ \beta (A)  \mu(A)
     +(q^*-1-\beta(A))  \mu^{(2)} (A)  \right].  \end{eqnarray*} Thus we obtain the same inequality as (\ref{quadratic}) with $p$ replaced by $q^*. $ Following
     (\ref{quadratic}),  and repeating the same proof therein,
     we deduce that
     $$ q^* \geq 1 + \frac{
2\left[1-\alpha(A) \beta(A) \widetilde{\mu}(A)^2 \right]}
 {\mu^{(2)}(A) \left\{ \widetilde{\mu} (A)(\alpha(A)+\beta(A))
    + \sqrt{
 \left[\widetilde{\mu}(A) (\alpha(A)-\beta(A))
 \right]^2 +4}  \right\} },$$
 where $\widetilde{\mu} (A) : =  \mu (A)- \mu^{(2)} (A). $ ~~ $\Box.$

 \vskip 0.08in

 It is also worth briefly comparing the
Babel-function-based bound (\ref{Tropp}) and those developed in
section 2 of this paper. At a first glance, it seems that
(\ref{Tropp}) is more sophisticated than those developed in section
2. However,  two types of bounds are mutually independent in the
sense that one cannot definitely dominate the other in general. For
example, when $\alpha(A) \leq \widehat{q}-1 $ and $\alpha(A)
<1/\mu(A)$ where $\widehat{q}$ is defined in Proposition 3.2, we
have
$$ 1\leq \mu_1(\widehat{q}-1) =\max_{1\leq k\leq m} \sum_{j=2}^{\widehat{q}} |\widehat{G}_{kj}| \leq \alpha(A)
\mu(A)+(\widehat{q}-1-\alpha(A))\mu^{(2)}(A).$$ Thus,
$$\widehat{q} \geq 1+ \frac{1-\alpha(A) \widetilde{\mu}(A)}{\mu^{(2)}(A)} = \left( 1+ \frac{1}{\mu(A)}\right) +
  \left(\frac{1}{\mu^{(2)}(A)} -\frac{1}{\mu(A)} \right) (1-\alpha(A)
  \mu(A)). $$ In this case, the Babel-function-based bound  (\ref{Tropp}) is
  tighter than bound (\ref{unique-b}). However, when
  $ \widehat{q}-1< \alpha(A) $, the relationship between the bounds  (\ref{unique-b}) and (\ref{Tropp})
 can be complicated. The bound (\ref{unique-b}) and the one in Theorem 2.9 might be
   tighter than (\ref{Tropp}).    Indeed, let us assume that $ \widehat{q}-1< \alpha(A) \leq
  p-1$ where $p=\textrm{Spark}(A), $  and  $\alpha(A)< 1/\mu(A).$
  Then (\ref{unique-b})  indicates that
  $$p = \left \lceil 1+ \frac{1}{\mu(A)} +
  \left(\frac{1}{\mu^{(2)}(A)} -\frac{1}{\mu(A)} \right) (1-\alpha(A)
  \mu(A))\right\rceil + t^*$$ for some integer $t^*\geq 0.$ This can be
  written as
  $$ \widehat{q} = \left \lceil 1+ \frac{1}{\mu(A)} +
  \left(\frac{1}{\mu^{(2)}(A)} -\frac{1}{\mu(A)} \right) (1-\alpha(A)
  \mu(A))\right\rceil + t^*- (p-\widehat{q}).$$
If $t^* < p-\widehat{q}, $ then  the above inequality  implies that
$$ \widehat{q} \leq  1+ \frac{1}{\mu(A)} +
  \left(\frac{1}{\mu^{(2)}(A)} -\frac{1}{\mu(A)} \right) (1-\alpha(A)
  \mu(A)). $$
By Proposition 2.6, the right-hand side of the above is dominated by
$\Psi (\alpha(A), \beta(A), \mu(A), \mu^{(2)}(A)) $. Therefore, as a
lower bound of spark,   (\ref{lowerbound}) is tighter than
(\ref{Tropp}) in this case.

\section{Scaled mutual coherence }

The  mutual coherence is not only an important property for the
development of the uniqueness of sparsest solutions, but also
crucial for the performance guarantee and stability analysis for
many sparsity-seeking algorithms, such as basis pursuit, orthogonal
matching pursuit,  and thresholding algorithms (see e.g.,
\cite{CDS98, EB02, DE03, FN03, T04, T04b, DET06,  BDE09,  E10,
EK12}).    The Babel function \cite{DE03,T04}, fusion coherence
\cite{BRK11}  and block coherence \cite{EKB10} are several variants
of the mutual coherence. In this section, we introduce the scaled
mutual coherence, which may lead to an optimal coherence-based
estimate of the spark in certain sense. In theory, the improved
results established in previous sections can  be either extended or
further improved by choosing a suitable scaling matrix.

\subsection{Uniqueness via the scaled mutual coherence}
Note that Spark(A), where $A\in R^{m\times n} $ with $m<n,$  is
invariant under a nonsingular linear transformation in the sense
that
$$ \textrm{Spark}(A) = \textrm{Spark}(WA)$$ for any nonsingular
matrix $W\in R^{m\times m}. $  However, the mutual coherence $\mu(A)
$ is not. That is,
$$\mu (A) \not=\mu(WA)$$ in general (see Examples 4.4 and 4.5 in this section). Thus
 the improved  conditions    (\ref{unique-a}) -(\ref{Coherence-rank-1})
still have a room for a further improvement by using a suitable
nonsingular scaling $W.$ Motivated by this observation, we consider
the weighted inner product between every pair of columns of a
matrix, and define
$$ \mu_W (A) =\max _{i\not= j} \frac{|(Wa_i)^TW a_j|}{\|W a_i\|_2
\cdot \|W a_j\|_2} = \mu(WA). $$ Similarly, we define
$$\mu^{(2)}_W(A) = \max _{i\not= j}  \left\{ \frac{ |(Wa_i)^TW
a_j|}{\|W a_i\|_2 \cdot \|W a_j\|_2}: ~ \frac{|(Wa_i)^TW a_j|}{\|W
a_i\|_2 \cdot \|W a_j\|_2} < \mu_W(A) \right\} =\mu^{(2)}(WA).
$$ In this paper, $\mu_W(A) $ and $\mu^{(2)}_W(A) $ are referred to be as the scaled mutual
coherence and  the scaled sub-mutual coherence, respectively. It
makes sense to introduce the next definition.

\vskip 0.08in

\textbf{Definition 4.1.} \emph{Let $$ \mu_*(A) := \min_{W}
\left\{\mu_W(A): ~W \in R^{m\times m} \textrm{ is
nonsingular}\right\}.  $$ $\mu_*(A) $ is called the optimal scaled
mutual coherence (OSMC) of $A.$ }

 \vskip 0.08in

By definition, we have
 $ \mu_*(A) \leq \mu_W(A)$  for any nonsingular $W\in R^{m\times m} $ and any $A\in R^{m\times n}.$ In particular,
 by setting $W=I$ (the identity matrix),  we see that
 $\mu_*(A) \leq \mu(A)$ for any $A. $  As shown by the
 next result, the OSMC provides a theoretical lower bound
 for the spark that is better than any other scaled-mutual-coherence-based
bound.

 \vskip 0.08in

\textbf{Theorem 4.2.} \emph{For any $m\times n$ ($m<n$)  matrix $A$
with nonzero columns, we have  $\mu_*(A) >0, $ and
 $$   1+\frac{1}{\mu_W(A)} \leq  1+\frac{1}{\mu_*(A)} \leq   \textrm{Spark}
(A) $$  for any nonsingular matrix $W\in R^{m\times m}.$ Hence if
the system $Ax=b$ has a solution satisfying
$$\|x\|_0\leq   \left(1+\frac{1}{\mu_*(A)} \right)/2, $$   or more
restrictively, if there is a nonsingular matrix $W$ such that
$\|x\|_0<   \left(1+ 1/\mu_W(A) \right)/2, $  then $x$ is the unique
sparsest solution to the linear system.}

\vskip 0.08in

\emph{Proof.} Let $W$ be an arbitrary nonsingular matrix. We
consider the scaled matrix  $ WA.  $ Let $D=\textrm{diag}
(1/\|Wa_1\|_2, ..., 1/\|Wa_n\|_2) $ where $a_i, i=1,..., n$ are the
columns of $A.$  Then $WAD$ is a matrix with normalized columns.
Clearly, this normalization does not change the spark (and hence,
$\textrm{spark}(WAD)=\textrm{spark}(WA) =\textrm{spark}(A)$.) We
also note that $WAD (D^{-1} x) =Wb$ and $Ax=b$ have the same
sparsity of solutions. So without loss of generality, we assume that
all columns of $WA$ have unit $\ell_2$-norms.
 Let  $p=\textrm{Spark}(A).$ By definition, there exist $p$ columns of $A$
that are linearly dependent. Let $A_S$ consist of these $p$ columns.
  Then the matrix $$ G^{(W)}_{SS}:=  (WA_S)^T (WA_S)   = A_S^T W^TWA_S  $$   is
a $p\times p$ singular matrix due to the linear dependence of
columns of $A_S.$ Since $WA$ is normalized,  all diagonal entries of
$G^{(W)}_{SS}$ are equal to 1, and off-diagonal entries are less
than or equal to $\mu_W(A).$ By the singularity of $G^{(W)}_{SS}$,
this matrix has a zero eigenvalue. Thus by Gerschgorin's theorem,
there exists a row  of the matrix, say the $i$th row, such that
$$ 1 \leq \sum_{j\not=i} |(G^{(W)}_{SS})_{ij}| \leq (p-1) \mu_W(A), $$
which implied that
$$   \mu_W(A) \geq 1/(p-1) >0 .   $$
Note that the inequality above holds for any nonsingular matrix
$W\in R^{m\times m}.$ Taking the minimum value of the left-hand side
yields $\mu_*(A) \geq 1/(p-1) >0, $ and thus
$$ p (=\textrm{Spark}(A)) \geq 1+\frac{1}{\mu_*(A) }.  $$ The right-hand side of the inequality above
is greater than or equal $1+1/\mu_W(A)$ for any nonsingular $W,$
since $\mu_W(A) \leq \mu_*(A).$  The uniqueness of  sparsest
solutions of the linear system $Ax=b$ follows immediately from
Theorem 1.1. ~~ $\Box.$

\vskip 0.08in

 With a scaling matrix $W,$ we denote the scaled
coherence rank and
 scaled sub-coherence rank by
$\alpha_W(A) = \alpha (WA) $ and  $ \beta_W (A) = \beta (WA),$
respectively.  By applying the same proof of Theorem 2.5 to the
scaled matrix $WA,$ the lower bound of spark, together with
uniqueness conditions for sparsest solutions in section 2, can   be
stated in terms of $\mu_W(A), \mu^{(2)}_W(A), \alpha_W (A) $ and $
\beta_W(A) . $ First, we define a class of matrices as follows.
\begin{eqnarray}\label {BM}  \widetilde{{\cal M}}  & = & \Big{\{} A\in R^{m\times n}:  \textrm{there is a nonsingular
  } W \in R^{m\times m} \textrm{ such that } \nonumber \\
    & & ~~\textrm{ either }
\alpha_W(A) \leq \frac{1}{\mu_W(A)}  \textrm{ and } \beta_W (A)
  < \alpha_W(A),
 \textrm{ or }     \alpha_W(A) < \frac{1}{\mu_W(A)} \Big{\} } .
  \end{eqnarray}

 We
now state the counterpart of Theorems 2.9 and 2.10, and Corollary
2.12 via the scaled coherence and the scaled coherence rank.

 \vskip 0.08in

\textbf{Theorem 4.3.}  (i) \emph{ Let $W\in R^{m\times m} $ be a
nonsingular matrix such that $A \in \widetilde{{\cal M}}, $ defined
by (\ref{BM}). If the system $Ax=b$ has a solution $x$ such that
$$  \|x\|_0 < \frac{1}{2} \left[1 + \frac{
2\left(1-\alpha_W(A) \beta_W(A) \widetilde{\mu}_W(A)^2 \right)}
 {\mu^{(2)}_W(A) \left\{ \widetilde{\mu}_W (A)(\alpha_W(A)+\beta_W(A))
    + \sqrt{\left[\widetilde{\mu}_W(A) (\alpha_W(A)-\beta_W(A))
 \right]^2 +4}  \right\} }\right], $$
 where $\widetilde{\mu}_W (A) :=  \mu_W (A)- \mu^{(2)}_W (A),$ then $x$ is
 the unique sparsest solution to the linear system. In particular, the conclusion is valid if the following  condition holds
$$\|x\|_0  <  \frac{1}{2} \left[ 1+\frac{1}{\mu_W(A)} +
  \left(\frac{1}{\mu^{(2)}_W(A)} -\frac{1}{\mu_W(A)} \right) (1-\alpha_W(A)
  \mu_W(A))
  \right] .$$ }

(ii) \emph{Suppose that  $W\in R^{m\times m} $ is  nonsingular  such
that $\mu_W(A) <1 $ and  $\alpha_W(A)=1. $  Then if the solution $x$
of $Ax=b$ obeys
\begin{equation}\label{unique-3}
\|x\|_0    <     \frac{1}{2} \left[ 1+\frac{1}{\mu_W(A)} +
  \left(\frac{1}{\mu^{(2)}_W(A)} -\frac{1}{\mu_W(A)} \right) (1-
  \mu_W(A))
  \right] ,
\end{equation}   $x$ is the unique sparsest solution to the linear system.}

\vskip 0.08in

The next example shows that ${\cal M} \subset \widetilde{{\cal M}},
$ i.e., $\widetilde{{\cal M}} $ is strictly larger than ${\cal M},$
and hence Theorem 4.3 covers a broader class of matrices than its
counterparts in section 2.3, and by a suitable scaling, the result
of Theorem 4.3  can further improve the results in section 2. In
fact, when $\alpha(A) \leq \frac {1}{\mu(A)}$ does not hold (in
which case $A\not\in {\cal M}$),  the scaled matrix $WA$ may satisfy
the condition $\alpha_W(A) < \frac{1}{\mu_W(A)} $ (so $WA\in {\cal
M},$ and hence $ A\in \widetilde{{\cal M}}$), as shown by the next
example.

\vskip 0.08in

\textbf{ Example 4.4.} Consider the matrix (\ref{Remark}) given in
Remark 2.3.  For this matrix,  $\mu(A)=0.9239$,
$\mu^{(2)}(A)=0.7644,$ $\alpha(A) =2,$   $\beta(A)=1,$ and the bound
(\ref{coherence}) is 1.0498.   Note that this matrix does not belong
to ${\cal M},$ since $\alpha(A) \not\leq 1/\mu(A).$ So
(\ref{unique-a})-(\ref{Coherence-rank-1}) cannot apply to this
matrix. Now, we randomly generate a scaling matrix as follows
$$ W= \left[
     \begin{array}{rrr}
       -0.9415 & -0.5320 & -0.4838 \\
       -0.1623 & 1.6821 & -0.7120 \\
       -0.1461 & -0.8757 & -1.1742 \\
     \end{array}
   \right].
     $$
It is easy to verify that $\mu_W(A)= 0.8954,$
$\mu^{(2)}_W(A)=0.8302,$ and $\alpha_W(A)=\beta_W(A) =1. $ In fact,
after this scaling, the absolute Gram matrix of the normalized $WA$
is given by
$$\textrm{abs}(G^{(W)})= \left[
  \begin{array}{cccccc}
  1.0000  &  0.3561  &  0.7138  &  0.8302   & 0.3978 &   0.8954\\
    0.3561  &  1.0000  &  0.5753  &  0.8130 &   0.7126 &   0.0973\\
    0.7138  &  0.5753  &  1.0000  &  0.8227  &  0.0177  &  0.4874\\
    0.8302  &  0.8130  &  0.8227  &  1.0000  &  0.1707  &  0.4969\\
    0.3978  &  0.7126  &  0.0177  &  0.1707  &  1.0000  &  0.7634\\
    0.8954  &  0.0973 &   0.4874  &  0.4969  &  0.7634  &  1.0000
  \end{array}
\right].$$ Thus by this scaling, the original coherence rank
$\alpha(A)=2$ is down to $\alpha_W(A)=1,$ and hence Theorems 4.3(ii)
can apply to $WA.$ Note that the scaled bound
$(1+\frac{1}{\mu_W(A)})/2$ in Theorem 4.2 is $1.0584$ and the bound
(\ref{unique-3}) in Theorem 4.3 (ii) is $1.0630. $ Both improve the
original unscaled  bound (\ref{coherence}). This example shows that
while $A\not\in {\cal M},$ we  have $ WA \in {\cal M}, $ and hence
$A \in \widetilde{{\cal M}}.$

\vskip 0.08in

From simulations, we observe that when the coherence rank of a
matrix is high in the sense that $\alpha(A) \geq 2,$ it is  quite
sensitive to a scaling $W,  $ which may immediately reduce
$\alpha(A) $ to $\alpha_W(A)=1,$ as shown by the above example. When
the coherence rank $\alpha(A) =1,$ it is insensitive to a scaling
$W,$ and it is highly likely that $\alpha_W(A)$ remains 1.

 \vskip 0.08in

 \textbf{Example 4.5.} Consider $A$ and the absolute Gram matrix
 $\textrm{abs}(G)$ of its normalized counterpart
 $$ A=\left[
      \begin{array}{rrrr}
        0.0010 & -0.7998 & -0.6002 & 1.4290 \\
        0.8001 & -0.3558 & 0.4798 & 1.2393 \\
        0.5999 & 0.4801 & -0.6398 & -0.6849\\
      \end{array}
    \right], ~ \textrm{abs}(G)=  \left[
                \begin{array}{cccc}
                  1 & 0.0025 & 0.0005 & 0.2894 \\
                  0.0025 & 1 & 0.0022 & 0.9523 \\
                  0.005 & 0.0022 & 1 & 0.0870 \\
                  0.2894 & 0.9523 & 0.0870 & 1 \\
                \end{array}
              \right].
 $$ For this example, $\alpha(A) =\beta(A) = 1, $ $ \mu(A) =0.9523,$
  and $
 \mu^{(2)}(A) = 0.2894.$ The standard  bound (\ref{coherence})
 is $(1+\frac{1}{\mu(A)})/2=1.025, $ which is improved to
 1.0824 by (\ref{Coherence-rank-1}). We now use the scaling matrix
 $$  W  = \left[
         \begin{array}{rrr}
           -0.2078 & 0.9393 & 0.1905 \\
           -0.9381 & 0.5715 & 0.3268   \\
           0.6702 & 0.2228 & 0.7662 \\
         \end{array}
       \right], $$
which is a randomly generated  nonsingular matrix. This scaling
matrix yields $\mu_W(A) =0.8343 ,$ $  \mu^{(2)}_A(A) =0.7272, $ and
$\alpha_W (A)=\beta_W(A) =1.$  The original  bound (\ref{coherence})
can be further improved by bound $(1+\frac{1}{\mu_W(A)})/2 =
1.0993,$ and by bound (\ref{unique-3}) that is equal to $1.1139.$ So
the scaled bound  improves the unscaled bound
(\ref{Coherence-rank-1}).

 \vskip 0.08in

The examples above do show that a scaling matrix can change the
mutual coherence, and may reduce coherence rank as well.  By a
suitable scaling, we can further improve the mutual-coherence-based
uniqueness conditions for sparsest solutions of  linear systems.
Note that if the OSMC is attainable, i.e., there exists a
nonsingular $W^*$ such that $\mu_*(A) =\mu_{W^*}(A).$ Then Theorem
4.3 holds for the OSMC. However, the optimal scaling $W^*$ is
difficult to obtain in general.  Also,   for a given linear system,
which scaling matrix  should be used in order to improve the
uniqueness claims is not obvious in advance.   However, the scaled
coherence can be viewed as a unified method for developing other
coherence-type conditions for the uniqueness of sparsest solutions.
It is worth mentioning that the Babel function can  be also
generalized to the weighted case, and related uniqueness claims can
be made as well.

\subsection{Application}

Note that the existing uniqueness claims for sparsest solutions of
linear systems are general and hold true uniformly for all $b.$
These claims are made largely by using the property of $A$ only, and
the role of $b,$ which is solution-dependent, has been overlooked.
Clearly, the property of the sparsest solution is usually dependent
on $A$ and $b.$  So it is interesting to incorporate the information
$b$ into a uniqueness criterion for sparsest solutions. The  scaled
mutual coherence can be used to achieve this goal. Indeed, let
$\phi$ be a mapping from $ R^m $ to $ R^m_{++} $ (the positive
orthant of $R^m$). Denote by $\Phi_u=\textrm{diag} (\phi(u)),$ a
nonsingular diagonal matrix with diagonal entries $\phi_i(u)>0 ,
i=1, ..., m.$ Setting $u=b,$ we see that  the system $Ax=b$ is
equivalent to
\begin{equation}\label{scaled-system} (\Phi_b A) x= \Phi_b b.
\end{equation}  For instance, we let $\phi(u) $ be separable, i.e.,
 $\phi(u) = (\phi_1(u_1), ..., \phi_n(x_n))^T,$ and  we define
\begin{equation} \label{phi} \phi_i(t) =
\left\{\begin{array}{cl}
1/t & \textrm{ if } t\not=0 \\
 1 & \textrm{otherwise}.
 \end{array} \right. \end{equation}
By this choice, we have $\Phi_b b = \textrm{diag} (\phi(b))b =
|\textrm{sign} (b)|.$
   Note that $\textrm{Spark}(A)
=\textrm{Spark} (\Phi_b A), $ and the sparsity of  solutions of the
scaled system (\ref{scaled-system})  is exactly the same as that of
$Ax=b.$ However, as we have seen before,  a  scaling matrix may
change the mutual coherence, and a suitable scaling may improve the
mutual-coherence-based uniqueness claims for sparsest solutions of a
linear system. Through a scaling matrix dependent on $b,$ the
contribution of $b$ to the uniqueness of sparsest solutions can be
demonstrated by the next two corollaries.

\vskip 0.08in

\textbf{Corollary 4.6.} \emph{
  If the system $ Ax=b,$ where $A\in R^{m\times n}$ with $m<n, $ has a solution satisfying
 $  \|x\|_0 < \left(1+\frac{1}{\mu (\Phi_b A)}\right)/2, $
then $x$ is the unique sparsest solution to the linear system.}

\vskip 0.08in

Applying to the scaled system (\ref{scaled-system}), this corollary
follows from Theorems 4.2 and 1.1 straightaway.  By Theorem 4.3,
this corollary can be  improved when the scaled coherence rank
$\alpha(\Phi_b A)$ is relatively small, as indicated  by the next
result.

 \vskip 0.08in

 \textbf{Corollary 4.7.}  \emph{Let $A$ be an $m\times n $ matrix  with $m<n.$}

(i) \emph{Suppose that either $ \alpha(\Phi_b A) \leq
\frac{1}{\mu(\Phi_b A)} $  and $ \beta(\Phi_b A)
  < \alpha(\Phi_b A) $ or
   $ \alpha (\Phi_b (A)) < \frac{1}{\mu(\Phi_b A)}.
  $     If the system $Ax=b$ has a solution $x$
  obeying  }
 \emph{{\small $$  \|x\|_0  <
\frac{1}{2} \left[1 + \frac{ 2\left(1-\alpha(\Phi_b A) \beta(\Phi_b
A) \widetilde{\mu}(\Phi_b A)^2 \right)}
 {\mu^{(2)}(\Phi_b A) \left\{ \widetilde{\mu}(\Phi_b A)(\alpha(\Phi_b A)+\beta(\Phi_b A))
    + \sqrt{
 \left[\widetilde{\mu}(\Phi_b A) (\alpha(\Phi_b A)-\beta(\Phi_b A))
 \right]^2 +4}  \right\} }\right], $$ }
 where $\widetilde{\mu}(\Phi_b A) :=  \mu(\Phi_b A)- \mu^{(2)}(\Phi_b
 A),$
then $x$ is the unique sparsest solution to the linear system. In
particular, the same conclusion holds if  $x$ obeys $$ \|x\|_0   <
\frac{1}{2} \left[ 1+\frac{1}{\mu(\Phi_b A)} +
  \left(\frac{1}{\mu^{(2)}(\Phi_b A)} -\frac{1}{\mu(\Phi_b A)} \right) (1-\alpha (\Phi_b A) \mu(\Phi_b
  A))
  \right],
$$
}

(ii)  \emph{If $\phi$ is chosen such that
     $\mu(\Phi_b A) <1 $ and
 $\alpha(\Phi_b A) =1,$   then the solution
$x$ of $Ax=b$ satisfying
\begin{equation}\label{unique-b3} \|x\|_0
< \frac{1}{2} \left[ 1+\frac{1}{\mu(\Phi_b A)} +
  \left(\frac{1}{\mu^{(2)}(\Phi_b A)} -\frac{1}{\mu(\Phi_b A)} \right) (1-
  \mu(\Phi_b A))
  \right]
\end{equation} is the unique sparsest solution of the linear system.}

\vskip 0.08in

 The next example
  shows that when $b$ is involved,  the uniqueness claim for sparsest
solutions can be improved in some situations.

\vskip 0.08in

\textbf{Example 4.8.} Consider the system $Ax =b$ where $A$ is a
$3\times 5$ matrix given by
$$ A=\left[
  \begin{array}{rrrrr}
    1 & -3 & -6 & 4 & -3 \\
    2 & 3 & -2 & -2 & 3 \\
    3 & -2 & 1 & 0 & 4 \\
  \end{array}
\right], \textrm{abs}(G)=  \left[
          \begin{array}{ccccc}
            1 & 0.1709 & 0.2922 & 0 & 0.6875 \\
            0.1709 & 1 & 0.3330 & 0.8581 & 0.3656 \\
            0.2922 & 0.3330 & 1 & 0.6984 & 0.4285 \\
            0 & 0.8581 & 0.6984 & 1 & 0.6903 \\
            0.6875 & 0.3656 & 0.4285 & 0.6903 & 1 \\
          \end{array}
        \right],$$
where $\textrm{abs}(G)$ is the absolute Gram matrix of the
normalized  $A. $ From $\textrm{abs}(G)$, we see that
$\mu(A)=0.8581,$ $\mu^{(2)}(A)= 0.6984, $ and $\alpha (A) =\beta(A)
=1.$ Thus  the standard bound (\ref{coherence}) is $1.0827$, which
is improved to $1.1016$ by (\ref{Coherence-rank-1}). In order to see
which  $b$ can further improve these bounds, let us randomly
generate a vector $b$, for instance, $b= (3.6159, -3.5189,
2.6954)^T. $ Let $\phi$ be given by (\ref{phi}). Then the absolute
Gram matrix of the scaled matrix $ \Phi_b A$ with normalized columns
is given by
$$ \textrm{abs}(G (\Phi_bA))=  \left[
          \begin{array}{ccccc}
             1.0000   & 0.3180  &  0.1608  &  0.0107  &  0.7833 \\
    0.3180  &  1.0000  &  0.2454  &  0.8042  &   0.1178\\
    0.1608   & 0.2454  &  1.0000  &  0.6784   & 0.4231 \\
    0.0107  &  0.8042  &  0.6784  &  1.0000  &  0.5928\\
    0.7833  &  0.1178  &  0.4231  &  0.5928  &  1.0000
          \end{array}
        \right],$$
from which we see that after this $b$-involved scaling,   the
coherence has changed to $\mu(\Phi_bA)=0.8042 $ and $\mu^{(2)}
(\Phi_bA)=0.7833, $ and
 the coherence rank remains unchanged. The scaled bound
$(1+\frac{1}{\mu(\Phi_b A)})/2 = 1.1217 $ and the scaled bound
(\ref{unique-b3}) equal to $1.1250 $ both improve the unscaled bound
(\ref{coherence}) and (\ref{Coherence-rank-1}).

 \vskip 0.08in

We now consider another  application of scaled mutual coherence.
Without loss of generality, we assume that $A$ is full-rank.  Let
$A=U \Sigma V^T$ be the singular value decomposition where $\Sigma$
is an $m\times m$ diagonal matrix with singular values as its
diagonal entries. We choose the scaling matrix $M= \Sigma^{-1} U^T$
which is nonsingular. Then $$ \textrm{Spark}(A) = \textrm{Spark}
(MA) = \textrm{Spark} (V^T) , ~ \mu(MA) =\mu(V^T).$$ This implies
that the lower bound of $\textrm{Spark}(A)$ can be computed by using
$V^T,$ instead of $A$ itself. Thus,  uniqueness claims for sparsest
solutions can be restated by the nonsquare orthogonal matrix $V^T.$
For completeness, we summarize this  result as follows.

\vskip 0.08in

\textbf{Corollary 4.9.}  \emph{Let $A$ be an $m\times n$ full-rank
matrix with $m<n,$ and  let $A=U\Sigma V^T$ be a singular value
decomposition where $\Sigma$ is diagonal of singular values .}

(i) \emph{If the system $ Ax=b$ has a solution satisfying
 $  \|x\|_0 < \left(1+\frac{1}{\mu (V^T)}\right)/2 $, then the
 system has a unique sparest solution.}

(ii) \emph{Suppose that either $ \alpha(V^T) \leq \frac{1}{\mu(V^T)}
$ and $ \beta(V^T)
  < \alpha(V^T), $ or
   $ \alpha (V^T) < \frac{1}{\mu(V^T)}
  .$    If the system $Ax=b$ has a solution $x$
  obeying
$$ \|x\|_0 < \frac{1}{2} \left(1 + \frac{ 2\left(1-\alpha(V^T)
\beta(V^T) \widetilde{\mu}(V^T)^2 \right)}
 {\mu^{(2)}(V^T) \left\{ \widetilde{\mu}(V^T)(\alpha(V^T)+\beta(V^T))
    + \sqrt{
 \left[\widetilde{\mu}(V^T) (\alpha(V^T)-\beta(V^T))
 \right]^2 +4}  \right\} }\right), $$
 where $\widetilde{\mu}(V^T) : =  \mu(V^T)- \mu^{(2)}(V^T),$ then $x$ is the unique sparsest solution to the linear system.
In particular, the  conclusion is valid  if the following condition
holds
$$ \|x\|_0 < \frac{1}{2} \left[ 1+\frac{1}{\mu(V^T)} +
  \left(\frac{1}{\mu^{(2)}(V^T)} -\frac{1}{\mu(V^T)} \right) (1-\alpha(V^T)
  \mu(V^T))
  \right].
$$}

 (iii)
 \emph{If $\mu(V^T) <1 $ and  $\alpha(V^T)=1,$  and if a solution $x$ of $Ax=b$
obeys $$ \|x\|_0    <    \frac{1}{2} \left[ 1+\frac{1}{\mu(V^T)} +
  \left(\frac{1}{\mu^{(2)}(V^T)} -\frac{1}{\mu(V^T)} \right) (1-
  \mu(V^T))
  \right],
$$ then $x$ is the unique sparsest solution of the linear system.}

\vskip 0.08in

\textbf{Remark 4.10.}  Consider the sparsest solution of the linear
system in matrix form
$$  {\cal A}\bullet x = \sum_{i=1}^N x_i A_i = B,$$
where $A_i\in R^{m\times q}, i=1,..., N ,$  $B\in   R^{m\times q}$
are given matrices, and $mq< N. $ The above system can be written as
a linear system in vector form, by using
$$ A= [vec(A_1), ..., vec(A_N)],  ~ b= vec (B),$$
where $vec(A_i)$ is a vector obtained by stacking the columns of
$A^T_i$ on top of another. Then $A$ is an $(mq)\times N$ matrix. We
also assume that $A$ is normalized in the sense that
$\|vec(A_i)\|_2=1 $ for $i=1, ..., N.$ To comply with the matrix
form, we define the Gram matrix of the linear operator ${\cal A}$
as
$$ G({\cal A})= \left[
\begin{array}{ccc}
 \textrm{tr}(A_1^TA_1)  & \cdots  &   \textrm{tr}(A_1^TA_N)\\
  \vdots &   & \vdots \\
  \textrm{tr}(A_N^TA_1)  &  \cdots &  \textrm{tr}(A_N^TA_N) \\
     \end{array}
       \right],
$$
and the mutual coherence as
$$ \mu ({\cal A}) =\max_{i\not =j} \frac{|\textrm{tr}(A_i^TA_j)| }{ \|A_i\|_F \cdot \|A_j\|_F
} , $$ where $\|\cdot\|_F$  denotes the Frobenius norm. Similarly,
we define $$ \mu ^{(2)}({\cal A}) =\max_{i\not =j} \left\{
\frac{|\textrm{tr}(A_i^TA_j)| }{ \|A_i\|_F \cdot\|A_j\|_F} : ~
\frac{|\textrm{tr}(A_i^TA_j)| }{ \|A_i\|_F\cdot \|A_j\|_F }<  \mu
({\cal A}) \right\} $$ as the second largest coherence, and $\alpha
({\cal A})$ is the maximum number of absolute entries   equal to
$\mu({\cal A})$ in a row of  $G({\cal A}).$ Then the results in
previous sections can be easily transformed to  the sparsest
solution of a linear system in matrix form.





\section{A further improvement  via support overlap}

Many uniqueness conditions for sparsest solutions of a linear system
were derived  from Theorem 1.1 by using the lower bound of
 $\textrm{Spark} (A).$ In this section, we point out that Theorem 1.1 itself might be
  improved in some situations by the   support
overlap of solutions of a linear system, leading to an enhanced
spark-type uniqueness condition. We use $\textrm{Supp}(x)$ to denote
the support of $x$, i.e., $\textrm{Supp}(x) =\{i: x_i\not = 0\}.$

\vskip 0.08in

\textbf{Definition 5.1.}      \emph{The support overlap $S^*$ of the
solution of $Ax=b$ is the index set
$$ S^*= \bigcap_{x\in {\cal Y}}
 \textrm{Supp} (x), $$ where
${\cal Y} =\{x: Ax=b\} $, the solution set of the linear system.}

\vskip 0.08in

 Clearly, $S^*$ might be empty if there is no common
index for the support of solutions. However, when some columns of
$A$ are crucial, and they must be used for the representation of $b,
$ the support overlap $S^*  $ is nonempty  for these cases.

\vskip 0.08in

\textbf{Theorem 5.2.}   \emph{Let $S^*$ be the support overlap of
the solution of the   system $Ax=b. $  If the   system has a
solution $x$ obeying
\begin{equation}\label{spark+S*} \|x \|_0    <  \frac{1}{2} (|S^*| +
\textrm{Spark} ( A )) , \end{equation}
 then $x$ is the unique sparsest
solution of the linear system.}

\vskip 0.08in

\emph{Proof.}  Let $x$ be a solution of the  system $Ax=b$ obeying
(\ref{spark+S*}). We now prove that it is the unique sparsest
solution of the linear system. We assume the contrary that the
linear system has a solution $y \not= x$ with $\|y\|_0 \leq
\|x\|_0.$ Since $ A(y-x) =0,$
 which implies that the columns $a_i, i\in \textrm{Supp}(y-x)$ of $A$ are linearly
 dependent, we have
 \begin{equation} \label{SSS}   \left\|y-
x\right\|_0  = |\textrm{Supp} (y-x) | \geq  \textrm{Spark} (A) .
\end{equation} Note that for any $u,v \in R^n$, the value of $
\|\textrm{diag} (u)v\|_0 $ is the number of $i$'s such that $u_iv_i
\not =0. $ So it is easy to see that
$$ S^*=  \min   \left\{ \left\|\textrm{diag}(x) u
\right\|_0 :  x, u \in {\cal Y}  \right\}. $$ Thus, for any $u,v\in
R^n, $ we have $$\|u- v\|_0 \leq \|u\|_0+\|v\|_0-\| \textrm{diag}
(u)v\|_0, $$ and hence
 \begin{eqnarray}  \label{Ineq}  \left\|y-
x \right\|_0   & \leq & \left\|y \right\|_0+ \left\|x \right\|_0 -
\left\|\textrm{diag}(x)
y \right\|_0 \nonumber \\
& \leq  &  2 \left\| x \right\|_0 -\left\|\textrm{diag}(x )
y  \right  \|_0  \nonumber \\
 & \leq &  2 \|x \|_0
-|S^*|,
\end{eqnarray}
where the first inequality follows from $ \|y \|_0 \leq \|x\|_0$ and
the second inequality follows from  the fact $
\left\|\textrm{diag}(x ) y  \right  \|_0 \geq |S^*| $ for any $x,y
\in {\cal Y}.$
 It follows from  (\ref{SSS}) and (\ref{Ineq}) that $  2\|x
  \|_0  -   |S^*| \geq    \textrm{Spark} (A), $
which contradicts with (\ref{spark+S*}). Thus $x$ is the unique
sparsest solution of the linear system. ~~  $\Box $

\vskip 0.08in

 As a result, all previous mutual-coherence-type
uniqueness criteria for sparsest solutions of a linear system can be
further improved when the value of $|S^*|$ or its lower bound   is
available. Taking Theorem 2.10 (iii) as an example, we have the
following result.

\vskip 0.08in

\textbf{ Corollary 5.3.}   \emph{Let $A\in R^{m\times n}, $ where
$m<n,$  be a matrix with $\mu(A) <1 $ and $\alpha(A) =1. $  Suppose
that $|S^*| \geq \gamma^* $ where $\gamma^*$ is known. Then if the
system $Ax=b$ has a solution $x$ obeying
\begin{equation} \label{rank-1} \|x\|_0    <    \frac{1}{2} \left[\gamma^*+ \left(1+\frac{1}{\mu(A)}\right) +
  \left(\frac{1}{\mu^{(2)}(A)} -\frac{1}{\mu(A)} \right) (1-
  \mu(A))
  \right]  ,
\end{equation}  $x$ is the unique sparsest solution of the linear  system. }

\vskip 0.08in

When the support overlap $S^*$ is nonempty, we have $|S^*| \geq 1.$
All the aforementioned mutual coherence type bounds for uniqueness
of sparsest solutions can be further improved by at least 0.5. Such
an improvement  can be crucial, as shown by the next example.

\vskip 0.08in

\textbf{Example 5.4.} Consider the system $Ax=b$ where
$$ A = \left[
         \begin{array}{ccccc}
           -1 & 0 & -4 & 2 & 4 \\
           0 & -1 & -1 & 1 & 2 \\
           0 & 0 & -1 & 0 & 0 \\
         \end{array}
       \right], ~~ b= \left[
                     \begin{array}{c}
                       2 \\
                       1/2 \\
                       1/2 \\
                     \end{array}
                   \right] $$
Clearly, the last two columns are linearly dependent. So
$\textrm{Spark}(A) =2,$ and Theorem 1.1 cannot confirm the
uniqueness of any sparsest solution. However, note that the third
column of $A$ is vital and must be used to represent $b.$ This means
that $x_3\not= 0 $ for any solution of the linear system. So,
$|S^*|\geq 1 =\gamma^*. $   Note that the solution $x^* =(0,0, 1/2,
0, 0)^T $ satisfies that
$$ \|x\|_0 =1 < 1.5=  (\gamma^* + \textrm{Spark}(A))/2 \leq
(|S^*|+\textrm{Spark} (A))/2. $$ By Theorem 5.2, $x^*$ is the unique
sparsest solution of the linear system. This example shows that by
incorporating the support overlap $S^*, $  the result of Theorem 1.1
can be remarkably improved when $ S^* \not= \emptyset.  $

\section{Uniqueness via range  property of $A^T$}

The exact recovery of all $k$-sparse vectors in $R^n$ by a single
matrix $A$ is called the uniform recovery.  To uniformly recover
sparse vectors, some matrix properties should be imposed on $A.$ The
restricted isometry property (RIP) \cite{CT05} and null space
property \cite{CDD09, Z08} are two well-known conditions for the
uniform recovery. Recently,    the so-called range space property
(RSP) of order $k$ was proposed in \cite{Z12a},  which can also
characterize the uniform recovery. All uniform recovering conditions
imply that the linear system $Ax=y: = Ax^0$ has a unique sparsest
solution. In fact,  these conditions have more capability than just
ensuring the uniqueness of sparsest solutions of a linear system.
For instance, they also guarantee that a linear system has a unique
least $\ell_1$-norm solution, leading to the  equivalence between
  $\ell_0$- and $\ell_1$-minimization
 problems,  which is fundamental for the development of
 compressed sensing theory.
In this section, we briefly discuss and develop certain  more
relaxed range properties of $A^T$ that guarantee the uniqueness of
sparsest solutions. Our first range property is defined as follows,
which was first introduced in \cite{ZL12} for a theoretical analysis
of reweighted $\ell_1$-methods for the sparse solution of a linear
system.

 \vskip 0.08in

\textbf{Definition 6.1} (Range  Property (I)).   ~ \emph{Let $A$ be
a full-rank $m\times n$ matrix with $m<n.$ Let $B$ be an
$(n-m)\times n ~$ matrix consisting of the basis of the null space
of $A.$ $B^T$ is said to satisfy a range space property (RSP) of
order $k$
 with a constant $\rho >0 $ if $$ \|
\xi_{\overline{J}}\|_1 \leq \rho \|\xi_J\|_1$$ for all $\xi\in {\cal
R}(B^T),$ the range space of $B^T, $  where $J\subseteq \{1, ...,
n\}$ with $|J|=k$ is the indices of $k$ smallest absolute components
of $\xi, $ and $\overline{J}=\{1,...,n\}\backslash J. $ }

\vskip 0.08in

Based on the above definition, we have the next result.

\vskip 0.08in

 \textbf{Theorem 6.2.}  \emph{ Let
$A\in R^{m\times n} $ and $ B\in R^{(n-m)\times n} $ be full-rank
matrices satisfying $A B^T =0,$ where $m<n. $  Suppose that  $B^T$
has a RSP of order $(n-k).$   Then the solution $x$ of the  system
$Ax=b$ obeying $\|x\|_0\leq k/2$ is the unique sparsest solution of
the linear system.}

\vskip 0.08in

 \emph{Proof.}  First, under the condition of the theorem, we have the
 following statement (see e.g., Proposition 3.6 in \cite{ZL12}):
  \emph{$B^T $ has the RSP of order $(n-k)$ with a   constant
  $\rho>0$ if and only if $A$ has the NSP of order $k$ with the same
  constant
$\tau= \rho. $} Therefore, by the definition of NSP of order $k,$ we
have $ \|\eta_\Lambda \|_1 \leq \tau
 \|\eta_{\overline{\Lambda}}\|_1 $ for all $\eta \in {\mathcal N} (A) $ and
all $\Lambda \subseteq\{1, 2, ..., n\}$ with $|\Lambda|\leq k, $
where $\overline{\Lambda}  = \{i: i\notin \Lambda\}.$ This implies
that the solution $x$ with $\|x\|_0 \leq k/2$ must be unique.  In
fact, we note that two $(k/2)$-sparse solutions $x$ and $y$
 satisfy  $A (x-y)=0,$ i.e., $x-y \in {\cal N}(A).$ Let $\Lambda=\textrm{Supp}(x-y).$
 Since $x-y$ is at most $k$-sparse, we have $|\Lambda| \leq k.$  By the NSP of order $k,$ we
have  $$\|x-y\|_1 = \|(x-y)_\Lambda\|_1 \leq \tau
\|(x-y)_{\overline{\Lambda}}\|_1 =0, $$  which implies that $x=y.$
Thus the $(k/2)$-sparse solution is the uniqueness sparsest
 solution of the linear system.  ~$\Box$

\vskip 0.08in

 The above theorem impose range  property  on the basis of
the null space of $A$, instead of on $A$ itself. We now impose a
range property on $A$ directly.

\vskip 0.08in

 \textbf{Definition 6.3} (Range Property (II)). \emph{There exists an integer $k$ such that
for any disjoint subsets $\Lambda_1, \Lambda_2$ of $\{1,..., n\}$
with $|\Lambda_1|+|\Lambda_2|=
 k$ and $|\Lambda_2|\leq 1,$ the range space $ {\mathcal R}(A^T) $ contains a vector $\eta$ satisfying
 $\eta_i =1 $ for all  $ i \in \Lambda_1,$  $ \eta_i =-1
 $  for all $ i\in \Lambda_2, $ and $ |\eta_i| <1 $ for $i\notin \Lambda_1 \bigcup \Lambda_2.$}

\vskip 0.08in

The above definition is a  relaxed version  of the range property
introduced in \cite{Z12a}. Under the above range property (II), we
can prove the following result.

\vskip 0.08in

\textbf{Theorem 6.4.}  \emph{Suppose that $A\in R^{m\times n} $ with
$m<n$  satisfies the range property (II).  Then if the system $Ax=b$
has a solution obeying $\|x\|_0 \leq k/2,$  $x$ is the unique
sparsest solution of the linear system. }

\vskip 0.08in

\emph{Proof.} Under the range property (II), we first  prove that
any $k$ columns of $A$ are linearly independent.  In fact, let
$\Lambda = \{\gamma_1, ..., \gamma_k\}$  be an arbitrary subset of
$\{1,..., n\}$ with $|\Lambda|=k. $  We now prove that the columns
of $A_\Lambda$ are linearly independent.
  It is sufficient to show that $z_\Lambda=0$ is the only
solution to the system $A_\Lambda z_\Lambda=0. $  In fact, let us
assume $A_\Lambda z_\Lambda=0. $ Then $z=(z_\Lambda,
z_{\overline{\Lambda}} =0) \in R^n $ is in  ${\cal N} (A). $
Consider the disjoint sets $\Lambda_1=\Lambda,$ and $\Lambda_2
=\emptyset . $ By the range property (II), there exists a vector
$\eta \in {\mathcal R}(A^T) $ with $\eta_i=1$ for all $i\in
\Lambda_1 =\Lambda.$ By the orthogonality of ${\cal N}(A)$ and
${\cal R}(A^T),$ we have $$0= z^T \eta= z_\Lambda^T
\eta_\Lambda+z_{\overline{\Lambda}}^T \eta_{\overline{\Lambda}}
=z_\Lambda^T \eta_\Lambda,$$ which is nothing but
\begin{equation} \label{zzss}  z_{\gamma_1}+ z_{\gamma_2} + \cdots  +z_{\gamma_k} =0.
\end{equation}  Now we consider an arbitrary
 pair of disjoint sets:  $$\Lambda_1= \Lambda \backslash \{\gamma_i\}, ~~
\Lambda_2= \{\gamma_i\}, $$  which satisfy that
$|\Lambda_1|+|\Lambda_2|=k$ and $ |\Lambda_2| \leq 1.$ By the range
property (II), there exists an $\eta\in {\mathcal R}(A^T) $ with
$\eta_{\gamma_j}=1$ for every $j\not = i$  and $\eta_{\gamma_i}=-1.
$ Again, it follows from $z^T\eta=0$   that
$$     ( z_{\gamma_1} + \cdots + z_{\gamma_{i-1}} + z_{\gamma_{i+1}}
\cdots + z_{\gamma_k}) - z_{\gamma_i}=0,  $$ which holds for every
$i$ with $1\leq i\leq k. $ Combining these relations and
(\ref{zzss}) implies that   $z_{\gamma_i}=0$ for all $i=1, ..., k, $
i.e., $ z_\Lambda=0.$
 So any $k$ columns
of $A$ are linearly independent. This implies that  $k<
\textrm{Spark}(A).$ The desired result follows immediately from
Theorem 1.1. ~~  $\Box. $

\section{Conclusions}     Through such concepts as  sub-mutual
  coherence, scaled mutual coherence,  coherence rank,  and sub-Babel
function,  we have developed several new and improved sufficient
conditions for a linear system to have a unique sparsest solution.
The key result established in this paper claims that when the
coherence rank of a matrix is low, the mutual-coherence-based lower
bound for the spark of a matrix can be improved.    We have also
demonstrated that the scaled mutual coherence, which  yields a
unified uniqueness claim, may further improve the unscaled
coherence-based uniqueness conditions if a suitable scaling matrix
is used. The scaled mutual coherence enables us to integrate the
right-hand-side vector $b$ of a linear system, and the orthogonal
matrix out of the singular value decomposition of $A$  into a
uniqueness criterion for the sparsest solution of a linear system.
Moreover, the support overlap of solutions and certain range
property of a matrix also play an important role in the uniqueness
of  sparsest solutions.

 \end{document}